\documentclass[10pt]{amsart}
\usepackage{amssymb}
\usepackage{comment}
\theoremstyle{plain}
\newtheorem{thm}{Theorem}
\newtheorem{question}{Question}
\theoremstyle{definition}
\newtheorem{prop}{Proposition}
\newtheorem{cor}[prop]{Corollary}
\newtheorem{lem}[prop]{Lemma}

\newtheorem{defn}[prop]{Definition}

\DeclareMathOperator{\dom}{dom}

\DeclareMathOperator{\Sk}{Sk}
\DeclareMathOperator{\otp}{otp}
\DeclareMathOperator{\tr}{Tr}

\DeclareMathOperator{\acc}{acc}

\DeclareMathOperator{\Ch}{Ch}

\DeclareMathOperator{\pr}{Pr}
\DeclareMathOperator{\id}{id}
\DeclareMathOperator{\ran}{ran}

\DeclareMathOperator{\cf}{cf}
\DeclareMathOperator{\nacc}{nacc}

\newcommand{\sk}{\vskip.05in}
\newcommand{\restr}{\upharpoonright}

\newcommand{\subs}{\subseteq}

\numberwithin{equation}{section}

\begin{document}
\title{A note on strong negative partition relations}
\author{Todd Eisworth}
\address{Department of Mathematics\\
         Ohio University\\
         Athens, OH 45701\\
         U.S.A}
\email{eisworth@math.ohiou.edu}
\keywords{successors of singular cardinals, scales, square brackets partition relations, club guessing}
\subjclass{03E02}
 \thanks{The author acknowledges support from NSF grant DMS 0506063.}
\date{\today}
\begin{abstract}
We analyze a natural function definable from a scale at a singular cardinal, and use it to obtain some
 strong negative square-brackets partition relations at successors of singular cardinals.
The proof of our main result makes use of club-guessing, and as a corollary we obtain a fairly easy proof
of a difficult result of Shelah connecting weak saturation of a certain club-guessing ideal with strong
failures of square-brackets partition relations. We then investigate the strength of weak saturation of such ideals
and obtain some results on stationary reflection.
\end{abstract}
\maketitle

\section{Introduction}

Recall that the {\em square-brackets partition relation} $\kappa\rightarrow[\lambda]^\mu_\theta$ of
Erd\H{o}s, Hajnal, and Rado \cite{ehr}  asserts that for every function $F:[\kappa]^\mu\rightarrow\theta$
 (where $[\kappa]^\mu$ denotes the subsets of $\kappa$ of cardinality $\mu$),
there is a set $H\subs\kappa$ of cardinality $\lambda$ such that
\begin{equation}
\ran(F\restr[H]^\mu)\neq\theta,
\end{equation}
that is, the function $F$ omits at least one value when we restrict
it to $[H]^\mu$.  The negations of square-brackets partition relations are particularly interesting,
as such combinatorial principles have many applications outside of pure set theory.

This paper is primarily concerned with the combinatorial statement
\begin{equation}
\label{eqn6}
\lambda\nrightarrow [\lambda]^2_\lambda,
\end{equation}
for $\lambda$ the successor of a singular cardinal.  The assertion (\ref{eqn6}) states that there exists
 a function $F:[\lambda]^2\rightarrow\lambda$ with the property that
 \begin{equation}
 \ran(F\restr[A]^2)=\lambda
 \end{equation}
 for any unbounded subset $A$ of $\lambda$.  Traditionally, more descriptive language is used when
 discussing (\ref{eqn6}) ---  $F$ is called a
 {\em coloring}, and (\ref{eqn6}) says that we can color the pairs of ordinals from $\lambda$ using $\lambda$
 colors in such a way that all colors appear in any unbounded subset of $\lambda$. Thus, Ramsey's Theorem
 fails for $\lambda$ in an incredibly spectacular way.

 The question of whether $\lambda\nrightarrow[\lambda]^2_\lambda$ necessarily holds for $\lambda$
 the successor of singular cardinal is still open.  Much research has been devoted to this question
 (particularly by Shelah \cite{Sh:355}, \cite{Sh:365}, \cite{413}, \cite{535}) and
the related question of whether  such a $\lambda$ can be a
J\'{o}nsson cardinal. This has resulted in a complex web of
conditions that tightly
constrains what a potential counterexample can look like, but still no proof that a counterexample cannot
exist has emerged.

In this paper, we show (assuming $\lambda=\mu^+$ for $\mu$ singular) that in many situations there is a natural coloring $c:[\lambda]^2\rightarrow\lambda$
with the property that $c$ takes on almost every color on every unbounded $A\subs\lambda$,
We have used two qualifying phrases in the previous sentence. The first -- ``in many situations'' -- we leave
vague for now, although our theorem is general enough to cover the case where $\mu$ is singular of uncountable
cofinality.  The second qualifying phrase -- ``almost every color'' -- means ``the set of omitted colors
is small in the sense that it lies in a certain ideal associated with Shelah's theory of guessing
clubs''.   The proof that the coloring has the required characteristics is a blending of techniques
due to Todor{\v{c}}evi{\'c} (namely, the method of {\em minimal walks} \cite{minimal,stevobook,stevochapter}) and techniques due to Shelah (combinatorics associated
with {\em scales} \cite{Sh:355, myhandbook}).

We will actually we working with negative square-brackets partition relations much stronger than those discussed
in the first two paragraphs of the paper.  In particular, we will be investigating instances of the following, which
is a specific case of a much more general property introduced studied by Shelah in several works~\cite{cardarith, 413, 535, 572}.

\begin{defn}
Let $\lambda=\mu^+$ for $\mu$ a singular cardinal, and suppose $\sigma\leq\lambda$.  We say $\pr_1(\lambda,\lambda,\sigma,\cf(\mu))$
holds if there is a function $f:[\lambda]^2\rightarrow\sigma$ such that whenever $\langle t_\alpha:\alpha<\lambda\rangle$
is a sequence of pairwise disjoint elements of $[\lambda]^{<\cf(\mu)}$, then for any $\epsilon<\sigma$ we can find $\alpha<\beta<\lambda$
such that
\begin{equation}
f(\zeta,\xi)=\epsilon\text{ for all $\zeta\in t_\alpha$ and $\xi\in t_\beta$.}
\end{equation}
\end{defn}

Note that $\pr_1(\lambda,\lambda,\sigma,\cf(\mu))$ implies $\lambda\nrightarrow[\sigma]^2_\lambda$, as we may
take the $t_\alpha$ to be singletons.

This work continues research begun in~\cite{535} and is further continued in~\cite{819}, where we (together
with Shelah) address the question of
what happens at successors of cardinals of countable cofinality.

\section{Preliminaries}

The background material needed for our results divides fairly neatly into four categories: club-guessing, minimal walks,
scales, and elementary submodels. We handle each of these topics in turn.

\medskip

\noindent{\sf Guessing clubs}

\smallskip

Shelah's theory of club-guessing has been a key ingredient in many of his theorems establishing negative
square-brackets partition relations. The foundations of the theory can be found in \cite{Sh:365}, while \cite{413}, \cite{535},
and \cite{572} provide glimpses of how useful the material can be in combinatorial set theory. We will be
concerned with a particular type of club-guessing that has proven exceedingly relevant in
this context, as well as a certain club-guessing ideal --- the ideal $\id_p(\bar{C},\bar{I})$ --- that has
heretofore not received as much attention in the literature as it perhaps deserves. This ideal has a fairly
 complicated definition, and so our initial discussion will be concerned with its description.

\begin{defn}
Let $S\subs\lambda$ be a stationary set of limit ordinals.  We say that a sequence
$\bar{C}=\langle C_\delta:\delta\in S\rangle$ is an
 $S$--club system if $C_\delta$ is closed and unbounded in $\delta$ for each $\delta\in S$.
  We can extend this notion to sets containing successor ordinals by requiring that for successor $\delta$, $C_\delta$
   is a closed subset of $\delta$ containing the predecessor of $\delta$.
     In the special case where $S$ is all of $\lambda$, we note that $\langle C_\alpha:\alpha<\lambda\rangle$ is a
      {\em $C$-sequence} in the sense of Todor{\v{c}}evi{\'c}.
\end{defn}
The above definition makes no demands on the order-type or cardinality of $C_\delta$ --- all
 that is required is that it be closed and unbounded in $\delta$. We note as well that although the terms ``$\lambda$-club
 system'' and ``$C$-sequence (on $\lambda$)'' denote exactly the same sorts of objects, we will preserve the distinct
 terminology because these objects are used for two entirely different reasons --- a sequence $\langle C_\delta:\delta<\lambda\rangle$
 used solely for club-guessing purposes will be referred using the former notation, while we will use the latter notation if we
 intend to use the object exclusively for minimal walks. In the author's opinion, the possibility of confusion is far outweighed
 by the advantage of staying consistent with notation in the extant literature.

\begin{defn}
\label{idealdefn}
Let $\bar{C}=\langle C_\delta:\delta\in S\rangle$ be an $S$-club system for $S$ a stationary subset of
some cardinal~$\lambda$, and suppose $\bar{I}=\langle I_\delta:\delta\in S\rangle$ is a sequence such that
$I_\delta$ is an ideal on $C_\delta$ for each $\delta\in S$. The ideal $\id_p(\bar{C},\bar{I})$
consists of all sets $A\subs\lambda$ such that for some closed unbounded $E\subs\lambda$,
\begin{equation}
\delta\in S\cap E\Longrightarrow E\cap A\cap C_\delta\in I_\delta.
\end{equation}
\end{defn}

Notice that $\id_p(\bar{C},\bar{I})$ is a proper ideal if and only if the sequence $(\bar{C},\bar{I})$
guesses clubs in the sense that for every closed unbounded $E\subs\lambda$, the set of $\delta\in S$
such that
\begin{equation}
E\cap C_\delta\notin I_\delta
\end{equation}
is stationary.  This is weaker than the usual notions of club-guessing prevalent in the literature, which
would require that $E$ contains almost all members of $C_\delta$, as measured by $I_\delta$.  Note as well
that a set $A$ is $\id_p(\bar{C},\bar{I})$-positive if and only if for every closed unbounded $E\subs \lambda$,
the set of $\delta\in S$ for which
\begin{equation}
\label{blah}
E\cap A\cap C_\delta\notin I_\delta
\end{equation}
is stationary.

We will be interested in $\id_p(\bar{C},\bar{I})$ for a particular choice of $\bar{I}$ which we define
after establishing a bit of notation.

\begin{defn}
Suppose $C$ is closed and unbounded in $\delta$.  Then
\begin{equation}
\nacc(C):=\{\alpha\in C:\sup(C\cap\alpha)<\alpha\},
\end{equation}
and $\acc(C):=C\setminus\nacc(C)$. (The notation ``$\acc$'' and ``$\nacc$'' comes from ``accumulation points''
and ``non-accumulation points''.)
\end{defn}

The notation of the preceding definition allows us to state the particular form of club-guessing
of concern to us in this paper.

\begin{defn}
\label{sgooddef}
Suppose $\lambda=\mu^+$ for $\mu$ a singular cardinal, and $S$ is a stationary subset of $\{\delta<\lambda:\cf(\delta)=\cf(\mu)\}$.
We say that $(\bar{C},\bar{I})$ is {\em an $S$-good pair} if the following conditions are satisfied:
\begin{enumerate}
\item $\bar{C}=\langle C_\delta:\delta\in S\rangle$ is an $S$-club system
\sk
\item $\sup\{|C_\delta|:\delta\in S\}<\mu$
\sk
\item $\bar{I}=\langle I_\delta:\delta\in S\rangle$
\sk
\item for $\delta\in S$, $I_\delta$ is the ideal on $C_\delta$ generated by sets of the form
\begin{equation}
\{\gamma\in C_\delta:\gamma\in\acc(C_\delta)\text{ or }\cf(\gamma)<\alpha\text{ or }\gamma<\beta\}
\end{equation}
for $\alpha<\mu$ and $\beta<\delta$.
\item for every closed unbounded $E\subs\lambda$, we have
\begin{equation}
\{\delta\in S: E\cap C_\delta\notin I_\delta\}\text{ is stationary.}
\end{equation}
\end{enumerate}
\end{defn}

Some remarks are in order here. First, the term ``$S$-good pair'' is not a standard one, nor do we intend it to be ---
it simply captures the particular instance of club-guessing relevant for the problem at hand.  Second, the ideals $I_\delta$
in the preceding definition are designated by ``$J^{b[\mu]}_{C_\delta}$'' in Shelah's research, and at the suggestion
of the referee we have gone with simpler notation.  Third, Shelah has shown that $S$-good pairs exist
for every stationary $S\subs \{\delta<\lambda:\cf(\delta)=\cf(\mu)\}$  in the case where $\mu$ has uncountable cardinality.  This result follows immediately from Claim~2.6 on page~127 of~\cite{cardarith} (see part (5) of his Remark~2.6A); there are a few small errors in the proof given there, so a complete (and simpler) proof appears as Theorem~2 of the forthcoming~\cite{819}. It is
still an open problem whether $S$-good pairs exist when $\mu$ has countable cofinality. Finally, with regard
to condition~(5), we remark that the statement ``$E\cap C_\delta\notin I_\delta$'' means that for every $\tau<\mu$ and
$\epsilon<\delta$, there is a $\gamma$ in $E\cap\nacc(C_\delta)$ such that $\gamma>\epsilon$ and $\cf(\gamma)>\tau$.

The ideal $\id_p(\bar{C},\bar{I})$ is a proper ideal whenever $(\bar{C},\bar{I})$ is an $S$-good pair
 because of condition~(5) in Definition~\ref{sgooddef}. In this situation, the equation~(\ref{blah})
  means that for every $\tau<\mu$ and $\epsilon<\delta$, there is a $\gamma\in\nacc(C_\delta)\cap E\cap A$
   such that $\gamma>\epsilon$ and $\cf(\gamma)>\tau$.

\medskip

\noindent{\sf Minimal Walks}

\smallskip

Suppose now that $\bar{e}=\langle e_\alpha:\alpha<\lambda\rangle$ is a $C$-sequence for some cardinal
$\lambda$. Following Todor{\v{c}}evi{\'c}, given $\alpha<\beta<\lambda$ the {\em minimal walk from $\beta$ to $\alpha$ along $\bar{e}$}
is defined to be the sequence $\beta=\beta_0>\dots>\beta_{n+1}=\alpha$ obtained by setting
\begin{equation}
\beta_{i+1}=\min(e_{\beta_i}\setminus\alpha).
\end{equation}
The use of ``$n+1$'' as the index of the last step is deliberate, as the ordinal $\beta_n$ (the penultimate step)
is quite important in our proof.  The {\em trace} of the walk from $\beta$ to $\alpha$ is defined by
\begin{equation}
\tr(\alpha,\beta)=\{\beta=\beta_0>\beta_1>\dots>\beta_n>\beta_{n+1}=\alpha\}.
\end{equation}

We make use of standard facts about minimal walks.  In particular, suppose $\delta$ is a limit ordinal, $\delta<\beta<\lambda$,
and  $\beta=\beta_0>\dots>\beta_{n+1}=\delta$ is the minimal walk from $\beta$ to $\delta$. For $i<n$,
we know that $\alpha\notin e_{\beta_i}$, and so
\begin{equation}
\gamma^*:=\max\{\max(e_{\beta_i}\cap\delta):i<n\}<\delta.
\end{equation}
Suppose now that $\gamma^*<\alpha<\delta$, and let $\beta=\beta_0^*>\dots>\beta^*_{n^*+1}=\alpha$ be the
minimal walk from $\beta$ to $\alpha$.  From the definition of $\gamma^*$ it follows that
\begin{equation}
\beta_i=\beta_i^*\text{ for }i\leq n.
\end{equation}
Thus, the walks from $\beta$ to $\delta$ and from $\beta$ to $\alpha$ agree up to and and including
 the step before the former reaches its destination.  A proper discussion of minimal walks and their applications
 is beyond the scope of this paper, and we develop the theory only to the degree that we need it
 for our proof. We refer the reader to \cite{minimal},~\cite{bekkali},\cite{stevochapter}, or~\cite{stevobook}
 for more information.

\medskip

\noindent{\sf Scales}

\medskip

The next ingredient we need for our theorem is the concept of a scale for a singular cardinal.

\begin{defn}
Let $\mu$ be a singular cardinal. A {\em scale for
$\mu$} is a pair $(\vec{\mu},\vec{f})$ satisfying
\begin{enumerate}
\item $\vec{\mu}=\langle\mu_i:i<\cf(\mu)\rangle$ is an increasing sequence of regular cardinals
such that $\sup_{i<\cf(\mu)}\mu_i=\mu$ and $\cf(\mu)<\mu_0$.
\item $\vec{f}=\langle f_\alpha:\alpha<\mu^+\rangle$ is a sequence of functions such that
\begin{enumerate}
\item $f_\alpha\in\prod_{i<\cf(\mu)}\mu_i$.
\item If $\gamma<\delta<\beta$ then $f_\gamma<^* f_\beta$, where  the notation $f<^* g$  means that $\{i<\cf(\mu): g(i)\leq f(i)\}$ is bounded in $\cf(\mu)$.
\item If $f\in\prod_{i<\cf(\mu)}\mu_i$ then there is an $\alpha<\beta$ such that $f<^* f_\alpha$.
\end{enumerate}
\end{enumerate}
\end{defn}

It is an important theorem of Shelah~\cite{Sh:355} that scales exist for any singular $\mu$; readers seeking
a gentler exposition of this and related topics should consult~\cite{hsw}, \cite{cummings}, or~\cite{myhandbook}.
If $\mu$ is singular and $(\vec{\mu},\vec{f})$ is a scale for $\mu$, then there is a
natural way to color the pairs of ordinals $\alpha<\beta<\mu^+$ using $\cf(\mu)$ colors, namely
\begin{equation}
\label{deltadef}
\Gamma(\alpha,\beta)=
\sup(\{i<\cf(\mu):f_\beta(i)\leq f_\alpha(i)\})
\end{equation}
Although we do not use it, we mention that the coloring $\Gamma$ defined above
witnesses that $\pr_1(\mu^+,\mu^+,\cf(\mu),\cf(\mu))$ holds for any singular $\mu$ --- this is
a result of Shelah (Conclusion~4.1(a) on page~67 of~\cite{Sh:355}), and Section~5 of~\cite{myhandbook} contains
an explication of this result.

The next lemma will be used in the proof of our main theorem; it is a result that
holds for arbitrary scales at singular cardinals~$\mu$.
 We remind the reader that notation of the form ``$(\exists^*\beta<\lambda)\psi(\beta)$''
  means $\{\beta<\lambda:\psi(\beta)\text{ holds}\}$ is unbounded below~$\lambda$, while
    ``$(\forall^*\beta<\lambda)\psi(\beta)$'' means that $\{\beta<\lambda:\psi(\beta)\text{ fails}\}$ is bounded
     below~$\lambda$.

\begin{lem}
\label{lem1}
Let $\lambda=\mu^+$ for $\mu$ a singular cardinal of cofinality $\kappa$, and suppose $(\vec{\mu},\vec{f})$ is a scale
for~$\mu$.  Then there is a closed unbounded $C\subs\lambda$ such that the following holds for every $\beta\in C$:
\begin{equation}
\label{eqn15}
(\forall^*i<\kappa)(\forall\eta<\mu_i)(\forall\nu<\mu_{i+1})(\exists^*\alpha<\beta)\left[f_\alpha(i)>\eta\wedge f_\alpha(i+1)>\nu\right].
\end{equation}
\end{lem}
\begin{proof}
Our first step is to prove the following statement:
\begin{equation}
\label{eqn14}
(\forall^*i<\kappa)(\forall\eta<\mu_i)(\forall\nu<\mu_{i+1})(\exists^*\alpha<\lambda)\left[f_\alpha(i)>\eta\wedge f_\alpha(i+1)>\nu\right].
\end{equation}
Assume by way of contradiction that~(\ref{eqn14}) fails.  It follows that there is an unbounded $I\subs\kappa$
(without loss of generality satisfying $i\in I\rightarrow i+1\notin I$) such that a ``bad pair'' $(\eta_i,\nu_i)$ exists
for every $i\in I$.  Define a function $f\in\prod_{i<\kappa}\mu_i$ by
\begin{equation}
f(i)=
\begin{cases}
\eta_i  &i\in I,\\
\nu_{i-1} &i-1\in I,\\
0    &\text{otherwise.}
\end{cases}
\end{equation}
For each $i\in I$, there is an $\alpha_i<\lambda$ such that
\begin{equation}
\alpha\geq\alpha_i\Longrightarrow \neg(\eta_i<f_\alpha(i)\text{ and }\nu_i<f_\alpha(i+1)).
\end{equation}
We now choose an $\alpha^*<\lambda$ greater than each such $\alpha_i$ and such that $f<^*f_{\alpha^*}$,
and the contradiction is immediate.

To see that there is a closed unbounded set of $\beta<\lambda$ satisfying (\ref{eqn15}), we apply (\ref{eqn14})
to choose an $i^*<\kappa$ so large that
\begin{equation}
(\forall\eta<\mu_i)(\forall\nu<\mu_{i+1})(\exists^*\alpha<\lambda)\left[f_\alpha(i)>\eta\wedge f_\alpha(i+1)>\nu\right].
\end{equation}
holds for all $i$ satisfying $i^*\leq i<\kappa$.

Given such an $i$, for each each pair $(\eta,\nu)\in\mu_i\times \mu_{i+1}$ we let $A^i(\eta,\nu)$
be the set of $\alpha<\lambda$ satisfying
\begin{equation}
\eta<f_\alpha(i)\text{ and }\nu<f_\alpha(i+1).
\end{equation}
There are only $\mu$ sets of the form $A^i(\eta,\nu)$, and therefore the set of $\beta<\lambda$ that are
limit points of all $A^i(\eta,\nu)$ simultaneously is closed and unbounded in $\lambda$ as required.
\end{proof}

\medskip

\noindent{\sf Elementary Submodels}

\medskip

Our conventions regarding elementary submodels are fairly standard --- we assume that $\chi$ is a sufficiently
large regular cardinal and let $\mathfrak{A}$ denote the structure $\langle H(\chi),\in, <_\chi\rangle$
where $H(\chi)$ is the collection of sets hereditarily of cardinality less than $\chi$, and $<_\chi$ is some
well-order of $H(\chi)$.  The use of $<_\chi$ means that our structure $\mathfrak{A}$ has definable Skolem functions,
and we obtain the set of {\em Skolem terms} for $\mathfrak{A}$ by closing the collection of Skolem functions
under composition.

\begin{defn}
Let $B\subs H(\chi)$. Then $\Sk_{\mathfrak{A}}(B)$ denotes the Skolem hull of $B$ in the structure $\mathfrak{A}$.
More precisely,
\begin{equation*}
\Sk_{\mathfrak{A}}(B)=\{t(b_0,\dots,b_n):t\text{ a Skolem term for $\mathfrak{A}$ and }b_0,\dots,b_n\in B\}.
\end{equation*}
\end{defn}

The set $\Sk_{\mathfrak{A}}(B)$ is an elementary substructure of $\mathfrak{A}$, and it is the smallest such structure
containing every element of $B$.

The following technical lemma due originally to Baumgartner~\cite{jb} will provide a key ingredient for our proof.

\begin{lem}
\label{newcharlem}
Assume that $M\prec\mathfrak{A}$ and let $\sigma\in M$ be a cardinal.  If we define $N=\Sk_{\mathfrak{A}}(M\cup\sigma)$
then for all regular cardinals $\tau\in M$ greater than $\sigma$, we have
\begin{equation*}
\sup(M\cap\tau)=\sup(N\cap\tau).
\end{equation*}
\end{lem}
\begin{proof}
Given an $\alpha\in N\cap\tau$, we must produce a $\beta\geq\alpha$ in $M\cap\tau$.  Since $\alpha$ is in $N$,
there is a Skolem term $t$ and parameters $\alpha_0,\dots,\alpha_i,\beta_0,\dots,\beta_j$ such that
\begin{equation*}
\alpha=t(\alpha_0,\dots,\alpha_i,\beta_0,\dots,\beta_j)
\end{equation*}
where each $\alpha_\ell$ is less than $\sigma$ and each $\beta_\ell$ is an element of $M\setminus\sigma$.

Now define a function $F$ with domain $[\sigma]^{i+1}$ by
\begin{equation*}
F(x_0,\dots, x_i)=
\begin{cases}
t(x_0,\dots, x_i,\beta_0,\dots,\beta_j)  &\text{if this is an ordinal less than $\tau$},\\
0 &\text{otherwise.}
\end{cases}
\end{equation*}
The function $F$ is an element of $M$, and so $\beta:=\sup(\ran(F))$ is in $M$ as well.  Since $\tau$ is
a regular cardinal, it is clear that $\alpha\leq\beta<\tau$ as required.
\end{proof}

As a corollary to the above, we can deduce an important fact about {\em characteristic functions} of models,
which we define next.

\begin{defn}
\label{chardef}
Let $\mu$ be a singular cardinal of cofinality~$\kappa$, and let $\vec{\mu}=\langle\mu_i:i<\kappa\rangle$ be an increasing
sequence of regular cardinals cofinal in $\mu$.  If $M$ is an elementary submodel of $\mathfrak{A}$
such that
\begin{itemize}
\item $|M|<\mu$,
\item $\langle \mu_i:i<\cf(\mu)\rangle\in M$, and
\item $\kappa+1\subs M$,
\end{itemize}
then the {\em characteristic function of $M$ on $\vec{\mu}$} (denoted $\Ch^{\vec{\mu}}_M$) is the function
with domain $\kappa$ defined by
\begin{equation*}
\Ch^{\vec{\mu}}_M(i):=
\begin{cases}
\sup(M\cap\mu_i) &\text{if $\sup(M\cap\mu_i)<\mu_i$,}\\
0  &\text{otherwise.}
\end{cases}
\end{equation*}
If $\vec{\mu}$ is clear from context, then we suppress reference to it in the notation.
\end{defn}

In the situation of Definition~\ref{chardef}, it is clear that $\Ch^{\vec{\mu}}_M$ is an element of the product
 $\prod_{i<\kappa}\mu_i$, and furthermore, $\Ch^{\vec{\mu}}_M(i)=\sup(M\cap\mu_i)$ for all
sufficiently large $i<\kappa$.  We can now see that the following corollary follows immediately from Lemma~\ref{newcharlem}.

\begin{cor}
\label{skolemhulllemma}
Let $\mu$, $\kappa$, $\vec{\mu}$, and $M$ be as in Definition~\ref{chardef}.
If $i^*<\kappa$ and we define $N$ to be $\Sk_{\mathfrak{A}}(M\cup\mu_{i^*})$,
then
\begin{equation}
\Ch_M\restr [i^*+1,\kappa)=\Ch_N\restr [i^*+1,\kappa).
\end{equation}
\end{cor}

We introduce one more bit of notation concerning elementary submodels, with an eye toward simplifying
the terminology used in various proofs throughout the paper.

\begin{defn}
Let $\lambda$ be a regular cardinal. A $\lambda$-approximating sequence is a continuous $\in$-chain
$\mathfrak{M}=\langle M_i:i<\lambda\rangle$ of elementary submodels of $\mathfrak{A}$ such that
\begin{enumerate}
\item $\lambda\in M_0$,
\item $|M_i|<\lambda$,
\item $\langle M_j:j\leq i\rangle\in M_{i+1}$, and
\item $M_i\cap\lambda$ is a proper initial segment of $\lambda$.
\end{enumerate}
If $x\in H(\chi)$, then we say that $\mathfrak{M}$ is a $\lambda$-approximating sequence over $x$ if
$x\in M_0$.
\end{defn}

Note that if $\mathfrak{M}$ is a $\lambda$-approximating sequence and $\lambda=\mu^+$, then $\mu+1\subs M_0$ because
of condition~(4) and the fact that $\mu$ is an element of each $M_i$.

\section{Preliminary results}

In this section, we will examine how minimal walks can be made to interact with $S$-good pairs.
The first thing we prove is that such $S$-good pairs can be ``swallowed'' by $C$-sequences --- this technique
is due to Shelah~\cite{cardarith} and it is used as well in~\cite{535}. The $C$-sequence
that results from such an operation possesses a weak form of coherence that is the key ingredient in our
proof.

\begin{lem}
\label{lem2}
Let $\lambda=\mu^+$ for $\mu$ a singular cardinal, and suppose $(\bar{C},\bar{I})$ is
an $S$-good pair for some stationary subset $S$ of $\{\delta<\lambda:\cf(\delta)=\cf(\mu)\}$.
There is a $C$-sequence $\bar{e}=\langle e_\alpha:\alpha<\lambda\rangle$  such that $|e_\alpha|<\mu$
for all $\alpha<\lambda$, and
\begin{equation}
\delta\in e_\alpha\cap S\Longrightarrow C_\delta\subs e_\alpha.
\end{equation}
\end{lem}
\begin{proof}
Let $\sigma<\mu$ be a regular cardinal distinct from $\cf(\mu)$, and let $\langle e_\alpha^*:\alpha<\lambda\rangle$
be a $C$-sequence such that $\otp(e^*_\alpha)=\cf(\alpha)$ for $\alpha$ a limit ordinal, and such that
$e^*_{\alpha+1}=\{\alpha\}$. For each $\alpha<\lambda$, we define a sequence $\langle e_\alpha[i]:i<\sigma\rangle$
as follows:
\begin{itemize}
\item $e_\alpha[0]=e^*_\alpha$
\sk
\item $e_\alpha[i+1]$ is the closure in $\alpha$ of $e_\alpha[i]\cup\bigcup\{C_\delta:\delta\in e_\alpha[i]\cap S\}$
\sk
\item for limit $i<\sigma$,  $e_\alpha[i]$ is the closure in $\alpha$ of $\bigcup_{j<i}e_\alpha[j]$.
\sk
\end{itemize}
Now we define $e_\alpha$ to be the closure in $\alpha$ of $\bigcup_{i<\sigma}e_\alpha[i]$; it is straightforward
to see that $\langle e_\alpha:\alpha<\lambda\rangle$ has all of the required properties (note that $|e_\alpha|<\mu$
because we have a uniform bound on the cardinalities of the $C_\delta$'s, and that taking closures in this context
does not increase cardinality).
\end{proof}

We sometimes refer to conclusion (2)  of the preceding lemma by saying that the $C$-sequence $\bar{e}$ {\em swallows the $S$-good pair}~$(\bar{C},\bar{I})$.
Note as well that the above lemma is much more general than stated --- the ideals $\bar{I}$ are irrelevant, as
the proof requires only that  $\sup\{|C_\delta|:\delta\in S\}$ is less than~$\mu$.

Our first result extracts a key property of minimal walks along $C$-sequences that swallow
a given good pair. The result is implicit in Shelah's work, but we isolate it here as it provides an explanation
for many of the theorems he obtains in the final section of~\cite{Sh:365} as well as our Theorem~\ref{mainthm} below.

\begin{defn}
Let $\bar{e}=\langle e_\alpha:\alpha<\lambda\rangle$ be a $C$-sequence on some cardinal $\lambda$. Given a
subset $t$ of $\lambda$, we say that a limit ordinal $\beta^*$ is $t$-ok if $\beta^*<\min(t)$ and there
exists an ordinal $\gamma^*<\beta^*$ (which we refer to as a {\em witness} for $\beta^*$ and $t$) such that
for any $\alpha$ in the interval $(\gamma^*,\beta^*)$ and any $\xi\in t$, the walk from $\xi$ down to $\alpha$
end extends the walk from $\xi$ to $\beta^*$. (So in particular, we have $\beta^*\in\tr(\alpha,\xi)$.)
\end{defn}

\begin{lem}
\label{oklem}
Suppose $\lambda=\mu^+$ for $\mu$ singular, and let
$\bar{e}=\langle e_\alpha:\alpha<\lambda\rangle$ be a $C$-sequence
that swallows an $S$-good pair $(\bar{C},\bar{I})$ for some
stationary $S\subs \{\delta<\lambda:\cf(\delta)=\cf(\mu)\}$. If
$\delta\in S$ and $t\in[\lambda]^{<\cf(\mu)}$ satisfies
$\delta<\min(t)$, then $I_\delta$-almost every element of $C_\delta$
is $t$-ok.
\end{lem}
\begin{proof}
For each $\xi\in t$, we let the walk from $\xi$ to $\delta$ along $\bar{e}$ consist of the ordinals
\begin{equation}
\xi=\beta_0^\xi>\dots>\beta^\xi_{n(\xi)>}\beta^\xi_{n(\xi)+1}=\delta,
\end{equation}
so $\tr(\delta,\xi)=\{\beta^\xi_i:i\leq n(\xi)+1\}$. Next, we define (for each $\xi\in t$)
\begin{equation}
\tau(\xi):=\left|e_{\beta^\xi_{n(\xi)}}\right|<\mu
\end{equation}
and
\begin{equation}
\epsilon(\xi):=\max\{\max(e_{\beta^\xi_i}\cap\delta):i<n(\xi)\}<\delta.
\end{equation}
Finally, define
\begin{equation}
\tau:=\cf(\mu)+\sup\{\tau(\xi):\xi\in t\}
\end{equation}
and
\begin{equation}
\epsilon:=\sup\{\epsilon(\xi):\xi\in t\}.
\end{equation}

Since $|t|<\cf(\mu)=\cf(\delta)$, it is clear that $\tau<\mu$ and $\epsilon<\delta$ and so it suffices
to establish that any $\beta^*\in\nacc(C_\delta)$ satisfying $\beta^*>\epsilon$ and $\cf(\beta^*)>\tau$ is $t$-ok.

Let such a $\beta^*$ be given; we will give a witness $\gamma^*$ for $\beta^*$ and $t$.  For each $\xi\in t$, we know
\begin{equation}
\delta\in e_{\beta^\xi_{n(\xi)}}
\end{equation}
and hence
\begin{equation}
C_\delta\subs e_{\beta^\xi_{n(\xi)}}
\end{equation}
because $\bar{e}$ swallows the pair $(\bar{C},\bar{I})$. In particular, this means that $\beta^*$
is in every such $e_{\beta^\xi_{n(\xi)}}$ and, furthermore, since $\cf(\beta^*)>\tau$ we know
that $\beta^*$ cannot be an accumulation point of this set. Thus,
\begin{equation}
\beta^*\in\nacc\left(e_{\beta^\xi_{n(\xi)}}\right)\text{ for all }\xi\in t.
\end{equation}
We now define
\begin{equation}
\gamma^*:=\epsilon + \sup\left\{\max(e_{\beta^\xi_{n(\xi)}}\cap\beta^*):\xi\in t\right\}.
\end{equation}
Since $\beta^*$ is a limit ordinal of cofinality greater than $|t|$, we know that $\gamma^*<\beta^*$.

Suppose now that $\gamma^*<\alpha<\beta^*$.
Given $\xi\in t$, since
\begin{equation}
\epsilon(\xi)<\alpha<\delta
\end{equation}
it follows (by a standard ``minimal walk argument'') that the walk from $\xi$ to $\alpha$ commences with the sequence
\begin{equation}
\xi=\beta_0^\xi>\dots>\beta^\xi_{n(\xi)},
\end{equation}
and our choice of $\gamma^*$ guarantees that
\begin{equation}
\beta^*=\min(e_{\beta^\xi_{n(\xi)}}\setminus\alpha).
\end{equation}
Thus, $\beta^*$ is the next step past $\beta^\xi_{n(\xi)}$ in the walk from $\xi$ to $\alpha$, and from this
we conclude that $\beta^*$ is $t$-ok.
\end{proof}

It will become clear in the course of the paper that the preceding lemma is the cornerstone
for the proofs we will give.   As an illustration of how it can be used, we will prove
an easy theorem connecting $\id_p(\bar{C},\bar{I})$ with minimal walks.
In order to state the result, we need the following (modification of a) definition from
 Todor{\v{c}}evi{\'c}'s~\cite{stevobook}.

\begin{defn}
Let $\bar{e}$ be a $C$-sequence on some cardinal $\lambda$.  The {\em trace filter of the $C$-sequence
$\bar{e}$} is defined to be the filter on $\lambda$ generated by sets of the form
\begin{equation}
\label{tracefilter}
\bigcup\{\tr(\alpha,\beta):\{\alpha,\beta\}\in [A]^2\}
\end{equation}
for $A$ an unbounded subset of $\lambda$.
\end{defn}

This is slightly at odds with the terminology of~\cite{stevobook}, as Todor{\v{c}}evi{\'c} defines the trace filter to be the {\em normal} filter on $\lambda$ generated by sets of the form~(\ref{tracefilter}).  However, the changing of notation here is done with his blessing.  We remark that the question of whether or not a given $C$-sequence has a proper trace filter is a delicate one --- it is not something that happens automatically. We send the reader to Section~8.2 of~\cite{stevobook} for more information on these matters.

\begin{thm}
\label{thm1}
Let $\lambda=\mu^+$ for $\mu$ singular, and assume $(\bar{C},\bar{I})$ is an $S$-good pair for some stationary
$S\subs\{\delta<\lambda:\cf(\delta)=\cf(\mu)\}$.  If $\bar{e}$ is a $C$-sequence such that
 \begin{itemize}
 \item $|e_\alpha|<\mu$ for all $\alpha<\lambda$, and
 \sk
 \item $\bar{e}$ swallows $(\bar{C},\bar{I})$,
 \sk
 \end{itemize}
then the trace filter of $\bar{e}$ is disjoint to the ideal $\id_p(\bar{C},\bar{I})$. In particular, the trace
filter of $\bar{e}$ is a proper uniform filter on $\lambda$ disjoint to the non-stationary ideal.
\end{thm}

\begin{proof}
Let $A$ be an unbounded subset of $\lambda$. We prove that
\begin{equation*}
\Gamma_A:=\bigcup\{\tr(\alpha,\beta):\{\alpha,\beta\}\in [A]^2\}
\end{equation*}
has measure one with respect to the ideal $\id_p(\bar{C},\bar{I})$.  To see this, we define $E$ to be the closed
unbounded subset of $\lambda$ consisting of ordinals $\alpha$ for which $\alpha=\sup(A\cap\alpha)$, and show
that
\begin{equation}
\delta\in S\cap E\Longrightarrow (E\cap C_\delta)\setminus\Gamma_A\in I_\delta.
\end{equation}
It suffices to show for each $\delta\in S\cap E$ with $E\cap C_\delta\notin I_\delta$ that
the set $\Gamma_A$ contains almost every member of $E\cap C_\delta$, as measured by the ideal $I_\delta$.

This follows almost immediately from Lemma~\ref{oklem}. To see why, simply fix any $\beta>\delta$.
By Lemma~\ref{oklem}, we know that $I_\delta$-almost every member of $C_\delta$ is $\{\beta\}$-ok. If
$\beta^*\in E\cap C_\delta$ is $\{\beta\}$-ok, then for all sufficiently large $\alpha\in A\cap\beta^*$
we have $\beta^*\in\tr(\alpha,\beta)$ and therefore $\Gamma_A$ contains $I_\delta$-almost all members of $E\cap C_\delta$.
\end{proof}

\section{The main theorem}

We move now to the main theorem of this paper. Throughout this section, we adopt the following list of
assumptions:

\sk

\noindent{\sf Assumptions}

\sk
\begin{itemize}
\item $\lambda=\mu^+$ for $\mu$ singular of cofinality~$\kappa$
\sk
\item $(\bar{C},\bar{I})$ is an $S$-good pair for some stationary $S\subs\{\delta<\lambda:\cf(\delta)=\kappa\}$
\sk
\item $\bar{e}=\langle e_\alpha:\alpha<\lambda\rangle$ is a $C$-sequence that swallows $(\bar{C},\bar{I})$
\sk
\item $(\vec{\mu},\vec{f})$ is a scale for $\mu$
\sk
\item $\mathfrak{A}=\langle H(\chi),\in, <_\chi\rangle$ for some sufficiently large regular cardinal~$\chi$.
\sk
\end{itemize}

Our next task is to define a certain function $c:[\lambda]^2\rightarrow\lambda$ using walks associated
with the $C$-sequence $\bar{e}$.

\begin{defn}
Given $\alpha<\beta<\lambda$, let $\beta=\beta_0>\beta_1>\dots>\beta_n>\beta_{n+1}=\alpha$ list $\tr(\alpha,\beta)$
in decreasing order.  The function $c:[\lambda]^2\rightarrow\lambda$ is defined by setting $c(\alpha,\beta)$
equal to $\beta_m$, where $m\leq n+1$ is the least number for which
\begin{equation}
\Gamma(\alpha,\beta_m)\neq\Gamma(\alpha,\beta).
\end{equation}
\end{defn}

The function $c$ can easily be described in English: to calculate the value of $c(\alpha,\beta)$, we first
compute $\Gamma(\alpha,\beta)$, and then walk along $\bar{e}$ until we reach an ordinal $\beta_m$ where $\Gamma(\alpha,\beta_m)$
is different from $\Gamma(\alpha,\beta)$.  This ordinal $\beta_m$ is the value of $c(\alpha,\beta)$.

\begin{thm}
\label{mainthm}
Suppose $\langle t_\alpha:\alpha<\lambda\rangle$ is a sequence of pairwise disjoint subsets of~$\lambda$,
each of cardinality less than $\kappa$.  Then for $\id_p(\bar{C},\bar{I})$-almost all $\beta^*<\lambda$
the following holds:
\begin{equation*}
(\exists^*\beta<\lambda)(\forall^*i<\kappa)(\exists^*\alpha<\beta^*)(\forall\zeta\in t_\alpha)(\forall\xi\in t_\beta)\left[c(\zeta,\xi)=\beta^*\wedge\Gamma(\zeta,\xi)=i\right].
\end{equation*}
\end{thm}
\begin{proof}
We first prove the theorem under the assumptions that   $\alpha<\min(t_\alpha)$ and $\sup(t_\alpha)<\min(t_\beta)$
whenever $\alpha<\beta<\lambda$. Given this, let $A$ be the set of all $\beta^*<\lambda$ for which the conclusion of the theorem fails, and assume by way of contradiction that $A$ is $\id_p(\bar{C},\bar{I})$-positive.

Let $\mathfrak{M}=\langle M_i:i<\lambda\rangle$ be a $\lambda$-approximating sequence
over all the objects mentioned so far.
We define $E$ to be the set of $\delta<\lambda$ for which $M_\delta\cap\lambda=\delta$.  Since
$E$ is a closed unbounded subset of $\lambda$ and $A\notin\id_p(\bar{C},\bar{I})$, we know there is a $\delta\in S\cap E$ with
\begin{equation}
\label{eqn16}
A\cap E\cap C_\delta\notin I_\delta.
\end{equation}
We have assumed $\delta<\min(t_\delta)$, so for each $\xi\in t_\delta$ let
\begin{equation}
\xi=\beta_0^\xi>\dots>\beta_{n(\xi)}^\xi>\beta^\xi_{n(\xi)+1}=\delta
\end{equation}
list $\tr(\delta,\xi)$ in decreasing order.

Taking (\ref{eqn16}) together with  Lemma~\ref{oklem}, we can find $\beta^*\in A\cap E\cap C_\delta$ and $\gamma^*<\beta^*$
with $\cf(\beta^*)>\kappa$ so that whenever $\gamma^*<\alpha<\beta^*$ and $\xi\in t_\delta$ the walk along $\bar{e}$ from $\xi$ to $\alpha$
commences with the sequence
\begin{equation}
\xi=\beta^\xi_0>\dots>\beta^\xi_{n(\xi)}>\beta^*.
\end{equation}
We will now prove that the following statement holds:
\begin{equation}
\label{4.5}
(\forall^*i<\kappa)(\exists^*\alpha<\beta^*)(\forall\zeta\in t_\alpha)(\forall\xi\in t_\delta)\left[c(\zeta,\xi)=\beta^*\wedge\Gamma(\zeta,\xi)=i\right].
\end{equation}
If we combine the above with the fact that $\sup(M_{\beta^*+1}\cap\lambda)\in M_\delta\cap\lambda=\delta$,
we find that $M_{\beta^*+1}\models \beta^*\notin A$, which is a contradiction.

For each $\alpha<\lambda$, we define a function $f^{\min}_\alpha$ as follows:
\begin{equation}
f^{\min}_\alpha(i)=\min\{f_\zeta(i):\zeta\in t_\alpha\}.
\end{equation}
Since $|t_\alpha|<\kappa$ for each $\alpha$, it follows that
\begin{equation}
(\forall^*i<\kappa)\left[f_\alpha^{\min}(i)=f_{\min(t_\alpha)}(i)\right].
\end{equation}
In particular, the sequence $\langle f^{\min}_\alpha:\alpha<\lambda\rangle$ is a scale, and we note
that this new scale is an element of $M_0$.

An appeal to Lemma~\ref{lem1} gives us a closed unbounded set in $M_0$ as there. In particular, this closed
unbounded set is also an element of $M_{\beta^*}$ and so $\beta^*=\sup(M_{\beta^*}\cap\lambda)$ is necessarily
a member of this closed unbounded set. We conclude
\begin{equation*}
(\forall^*i<\kappa)(\forall\eta<\mu_i)(\forall\nu<\mu_{i+1})(\exists^*\alpha<\beta^*)\left[f^{\min}_\alpha(i)>\eta\wedge f^{\min}_\alpha(i+1)>\nu\right].
\end{equation*}
This means we can choose $i_0<\kappa$ such that
\begin{equation*}
i_0\leq i<\kappa\Longrightarrow(\forall\eta<\mu_i)(\forall\nu<\mu_{i+1})(\exists^*\alpha<\beta^*)\left[f^{\min}_\alpha(i)>\eta\wedge f^{\min}_\alpha(i+1)>\nu\right].
\end{equation*}

The next piece of the proof makes use of Skolem hulls. Let us define
\begin{equation*}
x:=\{\lambda, \mu, \kappa, (\vec{\mu},\vec{f}), S, \bar{e},\langle t_\alpha:\alpha<\lambda\rangle,\beta^*\},
\end{equation*}
and define $M$ to be the Skolem hull in $\mathfrak{A}$ of $x$ together with all ordinals
less than or equal to $\kappa$, that is,
\begin{equation*}
M:=\Sk_{\mathfrak{A}}(x\cup \kappa+1).
\end{equation*}

Since $|M|=\kappa<\mu_0$, it follows that
\begin{equation}
\Ch_M(i)=\sup(M\cap\mu_i)\text{ for all }i<\kappa,
\end{equation}
where $\Ch_M$ is the characteristic function of $M$ from Definition~\ref{chardef}.
We note that $M$ can be computed by taking a Skolem hull in the model $M_{\beta^*+1}$ and
therefore $M\in M_\delta$.  In particular,
\begin{equation}
\Ch_M\in M_\delta\cap \prod_{i<\kappa}\mu_i
\end{equation}
and $\Ch_M<^* f_\alpha$ for some $\alpha\in M_\delta\cap\lambda=\delta$.

Thus, we can find $i_1<\kappa$ such that
\begin{equation}
\Ch_M\restr [i_1,\kappa)<f_{\beta^\xi_j}\restr [i_1,\kappa) \text{ for all $\xi\in t_\delta$ and $i\leq n(\xi)$.}
\end{equation}
Finally, let $i_2<\kappa$ be such that $\cf(\beta^*)<\mu_{i_2}$ and define
\begin{equation}
i^*=\max\{i_0,i_1, i_2\}.
\end{equation}
We claim that if $i^*\leq i<\kappa$, then
\begin{equation}
\label{goal}
(\exists^*\alpha<\beta^*)(\forall\xi\in t_\delta)(\forall\zeta\in t_\alpha)\left[c(\zeta,\xi)=\beta^*\wedge\Gamma(\zeta,\xi)=i\right].
\end{equation}

Fix such an $i$, and define
\begin{equation}
N=\Sk_{\mathfrak{A}}(M\cup\mu_i).
\end{equation}
By Corollary~\ref{skolemhulllemma}, it follows that
\begin{equation}
\Ch_N\restr [i+1,\kappa)=\Ch_M\restr [i+1,\kappa),
\end{equation}
and therefore by our choice of $i_1$, for all $\alpha\in N\cap\lambda$ we have
\begin{equation}
\label{eqn2}
f_\alpha\restr [i+1,\kappa)<f_{\beta^\xi_j}\restr[i+1,\kappa)\text{ for all $\xi\in t_\delta$ and $j\leq n(\xi)$.}
\end{equation}

We now define
\begin{equation*}
\eta^*:=\sup\{f_{\beta^\xi_j}(i):\xi\in t_\delta\text{ and }j\leq n(\xi)\}
\end{equation*}
and
\begin{equation*}
\nu^*=f_{\beta^*}(i+1).
\end{equation*}
Clearly $\eta^*<\mu_i$ and $\nu^*<\mu_{i+1}$, and both are elements of $N$.  By our choice of $i_2$,
we know $\cf(\beta^*)\subs N$ and therefore $N\cap\beta^*$ is unbounded in $\beta^*$.  Given this,
from our choice of $i_0$ it follows that
\begin{equation}
(\exists^*\alpha<\beta^*)\left[\alpha\in N\wedge f_\alpha^{\min}(i)>\eta^*\wedge f^{\min}_\alpha(i+1)>\nu^*\right].
\end{equation}

Suppose now that $\alpha\in N$ satisfies $\gamma^*<\alpha<\beta$,  $f_\alpha^{\min(i)}>\eta^*$, and $f^{\min}_{\alpha}(i+1)>\nu^*$.
From (\ref{eqn2}) and the choice of $\eta^*$, we conclude
\begin{equation}
\label{eqn3}
\Gamma(\zeta,\beta^\xi_j)=i\text{ for all $\zeta\in t_\alpha$, $\xi\in t_\delta$, and $j\leq n(\xi)$}.
\end{equation}
By our choice of $\nu^*$, we know
\begin{equation}
\label{eqn4}
\Gamma(\zeta,\beta^*)\geq i+1\text{ for all $\zeta\in t_\alpha$.}
\end{equation}
The conjunction of (\ref{eqn3}) and (\ref{eqn4}) establishes
\begin{equation}
(\forall\zeta\in t_\alpha)(\forall\xi\in t_\delta)\left[c(\zeta,\xi)=\beta^*\wedge \Gamma(\zeta,\xi)=i\right].
\end{equation}
Thus, (\ref{4.5}) holds and a contradiction arises because $\beta^*$ was chosen to be in the set~$A$.

We now fulfill the promise made at the start of the proof by handling the case of an arbitrary sequence $\langle t_\alpha:\alpha<\lambda\rangle$.  Given such a sequence, we define an increasing function $\iota:\lambda\rightarrow\lambda$ such that
\begin{itemize}
\item $\alpha<\min(t_{\iota(\alpha)})$ for all $\alpha<\lambda$, and
\sk
\item $\alpha<\beta<\lambda\Longrightarrow \max(t_\alpha)<\iota(\beta)$.
\sk
\end{itemize}
If we define $s_\alpha= t_{\iota(\alpha)}$, then our work applies to the sequence $\langle s_\alpha:\alpha<\lambda\rangle$. In particular, there is a set $B$ in the filter dual to $\id_p(\bar{C},\bar{I})$
so that the conclusion of our theorem (as it applies to $\langle s_\alpha:\alpha<\lambda\rangle$)
 holds for every $\beta^*\in B$.

 Let $B^*$ be the set of all $\beta^*\in B$ that are closed under the function~$\iota$. Since $\id_p(\bar{C},\bar{I})$
 extends the non-stationary ideal, it is clear that $B^*$ is in the filter dual to $\id_p(\bar{C},\bar{I})$ and routine
 to check that the conclusion of the theorem (as it applies to $\langle t_\alpha:\alpha<\lambda\rangle$) holds
 for all $\beta^*\in B^*$.
\end{proof}

We will shortly deduce an interesting theorem of Shelah as a corollary to our main result, but to do this
we need to fix some terminology.

\begin{defn}
Let $I$ be an ideal on some set $A$, and let $\sigma$ be a cardinal.
The ideal $I$ is weakly $\sigma$-saturated if $A$ cannot be partitioned into $\sigma$ disjoint $I$-positive
sets, i.e., there is no function $\pi:A\rightarrow\sigma$ such that
\begin{equation*}
\pi^{-1}(i)\notin I
\end{equation*}
for all $i<\sigma$.
\end{defn}
It is clear that any maximal ideal is weakly $2$-saturated, so weakly saturated ideals are not very difficult to
find. The rest of this paper will demonstrate that the question of ``how weakly saturated is $\id_p(\bar{C},\bar{I})$?'' is quite important. We begin with the following, which follows from the work in Section~4 of Shelah's~\cite{Sh:365}.

\begin{cor}
\label{cor}
Suppose $\lambda=\mu^+$ for $\mu$ singular, and $(\bar{C},\bar{I})$ is an $S$-good pair for some stationary subset $S$ of
$\{\delta<\lambda:\cf(\delta)=\cf(\mu)\}$.  If $\id_p(\bar{C},\bar{I})$ is not weakly $\sigma$-saturated, then $\pr_1(\lambda,\lambda,\sigma,\cf(\mu))$
holds; in particular, we have  $\lambda\nrightarrow[\sigma]^2_\lambda$.
\end{cor}
\begin{proof}
Let $\pi:\lambda\rightarrow\sigma$ partition $\lambda$ into disjoint $\id_p(\bar{C},\bar{I})$-positive sets, and
define $f:[\lambda]^2\rightarrow\sigma$ by
\begin{equation*}
f(\alpha,\beta)=\pi(c(\alpha,\beta)).
\end{equation*}
Suppose $\langle t_\alpha:\alpha<\lambda\rangle$ is a family of disjoint subsets of $\lambda$ each of cardinality
less than $\cf(\mu)$, and let $\epsilon<\sigma$ be given. By Theorem~\ref{mainthm}, we can find $\alpha<\beta^*<\beta$
such that $\pi(\beta^*)=\epsilon$ and
\begin{equation*}
c(\zeta,\xi)=\beta^*\text{ for all }\zeta\in t_\alpha\text{ and }\xi\in t_\beta.
\end{equation*}
It is clear that $f$ is constant with value $\epsilon$ when restricted to $t_\alpha\times t_\beta$.
\end{proof}
Shelah's original proof of the above is much more difficult, as he starts with a partition of $\lambda$
into $\id_p(\bar{C},\bar{I})$-positive sets and uses this as a parameter to define his coloring, whereas
we use a scale to get a single ``master coloring'' that can be used (in the sense of the proof of Corollary~\ref{cor}) in
conjunction with any such partition.

\section{From $\mu$ to $\mu^+$}

Let $\mu$ be a singular cardinal, and suppose $(\bar{C},\bar{I})$ is an $S$-good pair for some stationary
$S\subs\{\delta<\mu^+:\cf(\delta)=\cf(\mu)\}$.  The results of the previous section focused our
attention on the degree of weak saturation possessed by $\id_p(\bar{C},\bar{I})$.  In this section
we get an improvement of Corollary~\ref{cor}.

If we assume that $\id_p(\bar{C},\bar{I})$ is not weakly $\mu$-saturated, then Corollary~\ref{cor}
tells us
\begin{equation}
\label{eqn9}
\pr_1(\mu^+,\mu^+,\mu,\cf(\mu))
\end{equation}
holds, and we immediately obtain the relation
\begin{equation}
\label{eqn10}
\mu^+\nrightarrow[\mu]^2_{\mu^+}.
\end{equation}

Now an elementary argument allows us to ``upgrade'' the relation (\ref{eqn10}) to the case of $\mu^+$ colors,
that is, we can easily obtain the stronger result
\begin{equation}
\label{eqn11}
\mu^+\nrightarrow[\mu^+]^2_{\mu^+}
\end{equation}
from (\ref{eqn10}). It is natural to ask if we can also upgrade (\ref{eqn9}) to obtain
\begin{equation}
\label{eqn12}
\pr_1(\mu^+,\mu^+,\mu^+,\cf(\mu)).
\end{equation}
We do not know if (\ref{eqn12}) follows from (\ref{eqn9}) in general, but we have as consolation the following new
theorem that tells us that~(\ref{eqn12}) can be obtained from the same hypotheses we use to obtain~(\ref{eqn9}).

\begin{thm}
Suppose $\lambda=\mu^+$ for $\mu$ a singular cardinal, and suppose $(\bar{C},\bar{I})$ is an $S$-good
pair for some stationary $S\subs\{\delta<\lambda:\cf(\delta)=\cf(\mu)\}$. Then $\pr_1(\mu^+,\mu^+,\mu^+,\cf(\mu))$
holds if $\id_p(\bar{C},\bar{I})$ is not weakly $\mu$-saturated.
\end{thm}
\begin{proof}
For each $\alpha<\lambda$ we fix surjection $g_\alpha$ from $\mu$ onto $\alpha$, and let $\pi:\lambda\rightarrow\mu$ give
a partition of $\lambda$ into disjoint $\id_p(\bar{C},\bar{I})$-positive sets.  We also fix a function
$h:\cf(\mu)\rightarrow\omega$ such that $h^{-1}(\{n\})$ is unbounded in $\cf(\mu)$ for each $n<\omega$, and let
$\bar{e}$ be a $C$-sequence swallowing $(\bar{C},\bar{I})$.

Given $\alpha<\beta$, let $\beta=\beta_0>\dots>\beta_n>\beta_{n+1}=\alpha$ list $\tr(\alpha,\beta)$ (where
we walk along $\bar{e}$, just as in the proof of Theorem~\ref{mainthm}) in decreasing
order, and let $i^*$ denote $\Gamma(\alpha,\beta)$. We define $m(\alpha,\beta)$ to be the least $m\leq n+1$
 with
\begin{equation}
\Gamma(\alpha,\beta_{m})\neq i^*,
\end{equation}
so that in terms of the coloring from Section~4, we have
\begin{equation}
c(\alpha,\beta)=\beta_{m(\alpha,\beta)}.
\end{equation}
We also define
\begin{equation}
k(\alpha,\beta)=
\begin{cases}
m(\alpha,\beta)-h(\Gamma(\alpha,\beta))& \text{ if $h(\Gamma(\alpha,\beta))\leq m(\alpha,\beta)$,}\\
0& \text{otherwise}
\end{cases}
\end{equation}
Finally, we define
\begin{equation}
\label{cstardef}
c^*(\alpha,\beta)= g_{\beta_{k(\alpha,\beta)}}\bigl(\pi\circ c(\alpha,\beta)\bigr).
\end{equation}

A formula like (\ref{cstardef}) surely deserves some explanation, so we will describe the coloring we use
in English.  Given $\alpha<\beta$, we first compute $i^*=\Gamma(\alpha,\beta)$ and note that $h(i^*)$
is some natural number.  We then walk from $\beta$ to $\alpha$ until the first place where $\Gamma$ changes.
This isolates $\beta_{m(\alpha,\beta)}=c(\alpha,\beta)$, and $\pi(c(\alpha,\beta))$ records the piece
of the partition that contains the ordinal $\beta_{m(\alpha,\beta)}$.  Next, we turn around and retrace $h(i^*)$
steps of the walk from $\beta$ to $\beta_{m(\alpha,\beta)}$ (so we are walking up, not down).  This takes us to an ordinal
\begin{equation}
\beta_{k(\alpha,\beta)}>\beta_{m(\alpha,\beta)}.
\end{equation}
Now to compute the value of $c^*(\alpha,\beta)$, we take the bijection between $\mu$ and
$\beta_{k(\alpha,\beta)}$ and apply it to the ordinal $\pi(\beta_{m(\alpha,\beta)})$.

We now prove that the coloring $c^*$ has the properties required by $\pr_1(\lambda,\lambda,\lambda,\cf(\mu))$.
Let $\langle t_i:i<\lambda\rangle$ be a sequence of pairwise disjoint elements of $[\lambda]^{<\cf(\mu)}$.  Given $\epsilon<\lambda$, we need to find $\alpha<\beta$ such that
$c^*\restr t_\alpha\times t_\beta$ is constant with value $\epsilon$.

Let $\mathfrak{M}=\langle M_i:i<\lambda\rangle$ be a $\lambda$-approximating sequence
over all the parameters accumulated so far.
Let $E=\{\delta<\lambda:M_\delta\cap\lambda=\delta\}$, and choose $\delta\in S\cap E$ for which $E\cap C_\delta\notin I_\delta$.
We may (and will) assume $\delta<\min(t_\delta)$, so we can apply Lemma~\ref{oklem} and find $\beta^*\in E\cap C_\delta$
such that that $\beta^*$ is $t_\delta$-ok with $\gamma^*<\beta^*$ acting as a witness. Next, since $\epsilon\in M_{\beta^*}\cap\lambda=\beta^*$, we can find $\varsigma^*<\mu$
for which
\begin{equation}
g_{\beta^*}(\varsigma^*)=\epsilon.
\end{equation}

Since $\beta^*\in E$, both $\gamma^*$ and $\varsigma^*$ are in $M_{\beta^*}$, and (again using the fact
that $\beta^*\in E$) it follows that there must be a stationary set $T\subs\lambda$ and $g:T\rightarrow\lambda$
such that for all $\beta\in T$,
\begin{itemize}
\item $\beta<g(\beta)$,
\sk
\item $\beta$ is $t_{g(\beta)}$-ok, with $\gamma^*$ acting as witness, and
\sk
\item $g_{\beta}(\varsigma^*)=\epsilon$.
\sk
\end{itemize}
By passing to a stationary subset if necessary, we can assume
\begin{equation}
\alpha<\beta\text{ in }T\Longrightarrow \sup(t_{g(\alpha)})<\beta,
\end{equation}
so that if we define
\begin{equation}
t^*_{\beta}= \{\beta\}\cup t_{g(\beta)}\text{ for }\beta\in T,
\end{equation}
the resulting family of sets $\{t^*_\beta:\beta\in T\}$ is pairwise disjoint.

We know that $\pi^{-1}(\varsigma^*)$ is $\id_p(\bar{C},\bar{I})$-positive, and so an application
of Theorem~\ref{mainthm} to the family $\langle t_\beta^*:\beta\in T\rangle$ gives us a $\beta\in T$
and $\beta^*<\beta$ such that $\pi(\beta^*)=\varsigma^*$, and
\begin{equation}
\label{eqn17}
(\forall^*i<\cf(\mu))(\exists^*\alpha<\beta^*)(\forall\zeta\in t^*_\alpha)(\forall\xi\in t^*_\beta)[c(\zeta,\xi)=\beta^*\wedge\Gamma(\zeta,\xi)=i].
\end{equation}
Now let $n^*$ be the length of the walk from $\beta$ to $\beta^*$. Since $h^{-1}(n^*)$ is unbounded in~$\cf(\mu)$,
we can use (\ref{eqn17}) to select $i^*<\cf(\mu)$ and $\alpha\in T\cap\beta^*$ such that $h(i^*)=n^*$, $\gamma^*<\alpha$, and
\begin{equation}
\label{eqn18}
(\forall\zeta\in t^*_\alpha)(\forall\xi\in t^*_\beta)[c(\zeta,\xi)=\beta^*\wedge\Gamma(\zeta,\xi)=i^*].
\end{equation}
We finish by proving that $c^*(\zeta,\xi)=\epsilon$ whenever $\zeta\in t_{g(\alpha)}$ and
 $\xi\in t_{g(\beta)}$.

Suppose $\zeta\in t_{g(\alpha)}$ and $\xi\in t_{g(\beta)}$, and let
$\xi=\beta_0>\beta_1>\dots>\beta_n=\zeta$ list $\tr(\zeta,\xi)$ in decreasing order. Since $c(\zeta,\xi)=\beta^*$,
we know that $\beta^*=\beta_{m(\zeta,\xi)}$.  Now
\begin{equation*}
\gamma^*<\zeta<\beta^*<\beta<\xi,
\end{equation*}
and so the choice of $\gamma^*$ implies that the walk from $\xi$ to $\zeta$ must pass through $\beta$ before
proceeding on to $\beta^*$.  Since
$h(\Gamma(\zeta,\xi))$ is the length of the walk from $\beta$ to $\beta^*$, we know
\begin{equation*}
\beta_{k(\zeta,\xi)}=\beta.
\end{equation*}
Since $g_\beta(\varsigma^*)=\epsilon$ and $\pi(\beta^*)=\varsigma^*$, we conclude
\begin{equation*}
c^*(\zeta,\xi)=\epsilon,
\end{equation*}
as required.
\end{proof}

In the next section, we show that in our usual context, if an ideal of the form $\id_p(\bar{C},\bar{I})$ is weakly $\mu$-saturated for $\mu$ strong limit singular then every stationary subset of $\{\delta<\mu^+:\cf(\mu)\neq\cf(\delta)\}$ reflects.

\section{On the weak saturation of $\id_p(\bar{C},\bar{I})$}

In this section, we directly address the question of weak saturation of ideals of the form $\id_p(\bar{C},\bar{I})$.
Our results do not require the full strength of Definition~\ref{sgooddef}. In particular,
the requirement (2), although important for the arguments in the preceding two
theorems, is not a necessary ingredient in the proofs of this section.
With this in mind, we offer the following definition in the same spirit as Definition~\ref{sgooddef}.

\begin{defn}
\label{sfairdef}
Suppose $\lambda=\mu^+$ for $\mu$ a singular cardinal, and $S$ is a stationary subset of $\{\delta<\lambda:\cf(\delta)=\cf(\mu)\}$.
We say that $(\bar{C},\bar{I})$ is {\em an $S$-fair pair} if the following conditions are satisfied:
\begin{enumerate}
\item $\bar{C}=\langle C_\delta:\delta\in S\rangle$ is an $S$-club sequence
\item $\bar{I}=\langle I_\delta:\delta\in S\rangle$
\item for $\delta\in S$, $I_\delta$ is the ideal on $C_\delta$ generated by sets of the form
\begin{equation}
\{\gamma\in C_\delta:\gamma\in\acc(C_\delta)\text{ or }\cf(\gamma)<\alpha\text{ or }\gamma<\beta\}
\end{equation}
for $\alpha<\mu$ and $\beta<\delta$.
\item for every closed unbounded $E\subs\lambda$, we have
\begin{equation}
\{\delta\in S: E\cap C_\delta\notin I_\delta\}\text{ is stationary.}
\end{equation}
\end{enumerate}
\end{defn}

It is clear that any $S$-good pair is also $S$-fair. In contrast to $S$-good pairs, it is known that $S$-fair
 pairs exist for any stationary $S\subs\{\delta<\mu^+: \cf(\mu)=\cf(\delta)\}$ regardless of the cofinality of $\mu$ --- the countable cofinality case is handled by Claim~2.8 of page~131 in~\cite{cardarith}, while another proof (yielding more information) can be found in~\cite{819}.

Definition~\ref{idealdefn} still applies, so we have a proper ideal $\id_p(\bar{C},\bar{I})$ associated with every $S$-fair pair.  We note the following facts about this ideal:

\begin{prop}
\label{observation}
Let $\lambda=\mu^+$ for $\mu$ singular, and suppose $(\bar{C},\bar{I})$ is an $S$-fair pair for some
stationary subset $S$ of $\{\delta<\lambda:\cf(\delta)=\cf(\mu)\}$.
\begin{enumerate}
\item The ideal $\id_p(\bar{C},\bar{I})$ is closed under unions of fewer than $\cf(\mu)$ sets.
\item There is an increasing sequence $\langle A_i:i<\cf(\mu)\rangle$ of sets in $\id_p(\bar{C},\bar{I})$ such that
\begin{equation*}
\lambda=\bigcup_{i<\cf(\mu)}A_i.
\end{equation*}
(So in particular, this ideal can never be normal.)
\item If $\cf(\mu)<\sigma=\cf(\sigma)<\mu$, then the ideal $\id_p(\bar{C},\bar{I})$ is closed under
increasing unions of length $\sigma$.
\end{enumerate}
\end{prop}
\begin{proof}
The first and third statements follow easily from the fact that each of the ideals $I_\delta$
is closed under such unions (see Observation 3.2 on page 139 of~\cite{cardarith}).
To see the second statement, let
\begin{equation*}
A_i:=\{\delta<\lambda:\cf(\delta)\leq\mu_i\},
\end{equation*}
(note that successor ordinals land in $A_0$).
\end{proof}
The result (3) of the preceding proposition essentially says that the ideal $\id_p(\bar{C},\bar{I})$ is {\em $\sigma$-indecomposable},
a notion that has long history in the literature (Section~2 of \cite{myhandbook} gives many references for
this notion). More precisely, for regular $\sigma$ an ideal $I$ is $\sigma$-indecomposable if and only if it is closed
under increasing unions of length $\sigma$. Since we will never consider the case of $\sigma$-indecomposability
for singular $\sigma$, we will take the conclusion of (3) as our definition of $\sigma$-indecomposability when
we use this terminology in the sequel.

Let us now assume that $\lambda=\mu^+$ for $\mu$ singular, and $(\bar{C},\bar{I})$ is an $S$-fair pair for some stationary $S\subs\{\delta<\lambda:\cf(\delta)=\cf(\mu)\}$.
Let $\kappa$ denote the cofinality of $\mu$, and let $\langle \mu_i:i<\kappa\rangle$ be an increasing sequence of regular cardinals cofinal
in~$\mu$. There is a natural function $g:\lambda\rightarrow\kappa$ defined by
\begin{equation}
g(\delta):=\text{ least $i$ such that }\cf(\delta)<\mu_i,
\end{equation}
and the image of $\id_p(\bar{C},\bar{I})$ under this function,  defined by
\begin{equation}
A\in g(\id_p(\bar{C},\bar{I}))\Longleftrightarrow g^{-1}(A)\in\id_p(\bar{C},\bar{I}),
\end{equation}
is a proper $\kappa$-complete ideal on $\kappa$ because $\id_p(\bar{C},\bar{I})$ itself is a proper
 $\kappa$-complete ideal.

This simple observation already sheds considerable light on
questions of weak saturation for the ideal $\id_p(\bar{C},\bar{I})$,
for if $\id_p(\bar{C},\bar{I})$ is weakly $\sigma$-saturated for
some $\sigma$, then the same is true for the ideal
$g(\id_p(\bar{C},\bar{I}))$ on~$\kappa$. Thus, if
$\id_p(\bar{C},\bar{I})$ is a maximal ideal (something whose
consistency is still open) then the cofinality of $\mu$ must be
countable or a measurable cardinal.  If $\id_p(\bar{C},\bar{I})$ is
weakly $\cf(\mu)$-saturated, then $g(\id_p(\bar{C},\bar{I}))$ is
 a $\cf(\mu)$-complete $\cf(\mu)$-saturated ideal on $\cf(\mu)$, and a well-known
 argument of Ulam tells us that in this case $\cf(\mu)$ cannot be a successor
 cardinal.  It is open whether or not the statement considered in
 the preceding sentence can ever occur, but the following theorem
 and corollary show us that under mild cardinal arithmetic
 assumptions it cannot.

\begin{thm}
\label{cf(mu)}
Suppose $\lambda=\mu^+$ for $\mu$ singular, and $(\bar{C},\bar{I})$ is an $S$-fair pair for some stationary
$S\subs\{\delta<\lambda:\cf(\delta)=\cf(\mu)\}$.  If $2^{\cf(\mu)}<\mu$, then there is a function $g:\lambda\rightarrow\cf(\mu)$
such that $g(\id_p(\bar{C},\bar{I}))$ is the ideal of bounded subsets of $\cf(\mu)$.
\end{thm}
\begin{proof}
Let $\kappa$ denote the cofinality of~$\mu$, and
let $\langle \mu_i:i<\kappa\rangle$ be an increasing sequence of regular cardinals cofinal in~$\mu$.
If $f:\kappa\rightarrow\kappa$ is increasing, then we define the function $\Phi_f:\lambda\rightarrow\kappa$ by
setting $\Phi_f(\delta)$ equal to the least $i$ such that $\cf(\delta)<\mu_{f(i)}$.

Given such an $f$, let $I_f$ denote the ideal $\Phi_f(\id_p(\bar{C},\bar{I}))$ on $\kappa$. We note that $I_f$ contains
all the bounded subsets of $\kappa$, and we will prove
that there is an $f$ such that $I_f$ is exactly the ideal of bounded subsets of $\kappa$.

Assume by way of contradiction that no such function exists. Then for each increasing
 $f:\kappa\rightarrow\kappa$, there is an unbounded $A_f\subs\kappa$ such that
\begin{equation}
\label{condition}
\Phi^{-1}_f(A_f)\in\id_p(\bar{C},\bar{I}).
\end{equation}
It follows that there is a closed unbounded $E_f\subs\lambda$ such that
\begin{equation}
\label{condition2}
\delta\in S\cap E_f\Longrightarrow E_f\cap C_\delta\cap\Phi^{-1}_f(A_f)\in I_\delta.
\end{equation}
Let $E\subs\lambda$ be the intersection of all these sets $E_f$. It is clear that $E$ is closed and unbounded
in $\lambda$ because we assumed $2^\kappa$ to be less than $\mu$.

By applying definition of ``$S$-fair'', we choose $\delta\in S\cap E$ for which $E\cap C_\delta\notin I_\delta$, and let $\langle \delta_i:i<\kappa\rangle$
be an increasing and continuous sequence of ordinals cofinal in $\delta$.  We define a function
 $f^*:\kappa\rightarrow\kappa$ as follows:

 Begin by setting $f^*(0)=0$.  If $i<\kappa$ and $f^*(j)<\kappa$ has been defined for all $j<i$, then
 we define $\eta=\sup\{f^*(j):j<i\}$.  Since $E\cap C_\delta\notin I_\delta$, we know there is an $\epsilon\in E\cap \nacc(C_\delta)$
 for which $\delta_\eta\leq\epsilon$ and $\mu_\eta\leq\cf(\epsilon)$. Choose $f^*(i)<\kappa$
 so that $\epsilon<\delta_{f^*(i)}$ and $\cf(\epsilon)<\mu_{f^*(i)}$.

\medskip

It is clear that the function $f^*$ defined above is an increasing function from $\kappa$ to~$\kappa$, and
therefore $E_{f^*}$ exists and $E$ is a subset of $E_{f^*}$.
Because of this, an appeal to (\ref{condition2}) tells us that
\begin{equation}
\label{condition3}
E\cap C_\delta\cap\Phi^{-1}_{f^*}(A_{f^*})\in I_\delta.
\end{equation}
However, for each $i<\kappa$ the interval $[\delta_{f^*(i)},\delta_{f^*(i+1)})$ contains an ordinal
$\epsilon$ from $E\cap\nacc(C_\delta)$ with $\mu_{f^*(i)}\leq \cf(\epsilon)<\mu_{f^*(i+1)}$, and therefore
for {\em every} unbounded $A\subs\kappa$ we have
\begin{equation}
\label{condition4}
E\cap C_\delta\cap\Phi^{-1}_{f^*}(A)\notin I_\delta.
\end{equation}
The conjunction of (\ref{condition3}) and (\ref{condition4}) gives us the required contradiction.
\end{proof}

The preceding theorem yields the following corollary which
strengthens an unpublished result of Shelah.

\begin{cor}
If $\lambda=\mu^+$ for $\mu$ singular with $2^{\cf(\mu)}<\mu$ and $(\bar{C},\bar{I})$ is an $S$-fair
pair for some stationary $S\subs\{\delta<\lambda:\cf(\delta)=\cf(\mu)\}$, then $\lambda$ can be partitioned
into $\cf(\mu)$ disjoint $\id_p(\bar{C},\bar{I})$-positive sets.
\end{cor}
\begin{proof}
Fix $g:\lambda\rightarrow\cf(\mu)$ as in Theorem~\ref{cf(mu)}. If we partition $\cf(\mu)$ into disjoint
sets $\{ A_i:i<\cf(\mu)\}$ of cardinality $\cf(\mu)$, then the sets $\{ g^{-1}(A_i):i<\cf(\mu)\}$
give us the required partition of $\lambda$.
\end{proof}

The result of Shelah referred to above appeared in an unpublished
preliminary version of~\cite{535}, and established under the same
hypotheses the existence of pairs $(\bar{C},\bar{I})$ for which the
corresponding ideal fails to be $\cf(\mu)$-saturated.

Our next goal is to establish a connection between ideals of the
sort we have been considering and reflection of stationary sets. In
particular, we will prove the following theorem.

\begin{thm}
\label{reflthm}
Suppose $\lambda=\mu^+$ for $\mu$ a strong limit singular, and let
$(\bar{C},\bar{I})$ be an $S$-fair pair for some stationary $S\subs\{\delta<\lambda:\cf(\delta)=\cf(\mu)\}$.
If $\id_p(\bar{C},\bar{I})$ is weakly $\mu$-saturated,
then every stationary subset of $\{\delta<\lambda:\cf(\delta)\neq\cf(\mu)\}$ reflects.
\end{thm}

The proof of Theorem~\ref{reflthm} will use ideas from many places. One of the main ingredients is a
variant of the following result of Shelah from~\cite{413}.

\begin{thm}[Shelah~Claim~2.9 of \cite{413}]
\label{413thm}
Suppose $\lambda=\mu^+$ for singular strong limit $\mu$, and $(\bar{C},\bar{I})$ is an $S$-fair
pair for some stationary $S\subs\{\delta<\lambda:\cf(\delta)=\cf(\mu)\}$.
If $\id_p(\bar{C},\bar{I})$ is weakly $\theta^+$-saturated for some $\theta<\mu$,
then $|\mathcal{P}(\lambda)/\id_p(\bar{C},\bar{I})|\leq 2^\theta<\mu$.
\end{thm}

In the course of our proof of Theorem~\ref{reflthm} we will prove a mild generalization of
 Theorem~\ref{413thm}, but this seemed a good opportunity to make a few remarks on this result of Shelah.

The conclusion of Theorem~\ref{413thm} tells us that under the assumptions given,
 the ideal $\id_p(\bar{C},\bar{I})$ is ``almost'' maximal, in the sense that the reduced product $\mathcal{P}(\lambda)/\id_p(\bar{C},\bar{I})$ is small (if it
were maximal, the reduced product would have size 2).  If $\mu$ is of uncountable cofinality, then
the fact that $\id_p(\bar{C},\bar{I})$ is $\cf(\mu)$-complete gives this conclusion some added strength ---
there is a $\cf(\mu)$-complete filter on $\lambda$ that is close to being an ultrafilter in some sense.
This explains some rather cryptic remarks in the Analytical Guide to \cite{cardarith} --- Shelah says (page 462 of \cite{cardarith})
that in the situation of the above theorem,
\begin{quote}
{\footnotesize\bf
$\lambda$ is close to being ``$\kappa$-supercompact''}.
\end{quote}
It seems that what is meant here is that  $\cf(\mu)$ is close to being $\mu^+$-supercompact, in the sense
that there is a uniform $\cf(\mu)$-complete filter on $\mu^+$ that is close to being an ultrafilter.

We turn now to the generalization of Theorem~\ref{413thm} that we require for our proof of Theorem~\ref{reflthm}.

\begin{lem}
\label{cardlem}
Suppose $\lambda=\mu^+$ for $\mu$ a strong limit singular, and let $(\bar{C},\bar{I})$ be an $S$-fair pair for some stationary $S\subs\{\delta<\lambda:\cf(\delta)=\cf(\mu)\}$.
If $\id_p(\bar{C},\bar{I})$ is weakly $\mu$-saturated, then there is an $\id_p(\bar{C},\bar{I})$-positive
set $A\subs\lambda$ so that $$|\mathcal{P}(A)/\id_p(\bar{C},\bar{I})|<\mu.$$
\end{lem}
\begin{proof}
From previous work, we know that there is a partition $\langle A_i: i<\cf(\mu)\rangle$ of $\lambda$ into
disjoint $\id_p(\bar{C},\bar{I})$-positive sets.  Since $\id_p(\bar{C},\bar{I})$ is assumed to be weakly
$\mu$-saturated, it is clear that for some $i<\cf(\mu)$ and $\tau\in[\cf(\mu),\mu)$, the ideal
$\id_p(\bar{C},\bar{I})\cap\mathcal{P}(A_i)$
must be weakly $\tau$-saturated. Thus, we can fix an $\id_p(\bar{C},\bar{I})$-positive set $A$ and $\tau<\mu$
for which
\begin{equation}
\label{tausaturated}
\id_p(\bar{C},\bar{I})\cap\mathcal{P}(A)\text{ is weakly $\tau$-saturated.}
\end{equation}

Our argument now follows the line of attack used by Shelah in his proof of Theorem~\ref{413thm}. In particular,
we will prove the following statement:
\begin{equation}
\label{413goal}
\bigl|\mathcal{P}(A)/\id_p(\bar{C},\bar{I})\bigr|\leq 2^\tau.
\end{equation}

By way of contradiction, assume (\ref{413goal}) fails.  Let $M$ be an elementary submodel of $H(\chi)$
for some sufficiently large regular cardinal $\chi$ such that
\begin{itemize}
\item $S$, $A$, $(\bar{C},\bar{I})$, and $\id_p(\bar{C},\bar{I})$ are all in $M$,
\item $|M|=(2^\tau)^+$, and
\item $M$ is closed under sequences of length less than or equal to $\tau$.
\end{itemize}
It should be clear that such models exist.

Because $\mu$ is a strong limit cardinal, it follows that $|M|<\mu$, and therefore the set
\begin{equation}
E^*:=\bigcap\{E\in M:E\text{ closed unbounded in }\lambda\}
\end{equation}
is a closed unbounded subset of $\lambda$.  The set $E^*$ provides a litmus test for those subsets
of $\lambda$ that are in $M$, as a set $B\in M\cap\mathcal{P}(A)$ is $\id_p(\bar{C},\bar{I})$-positive if and only if there is
a $\delta\in S$ such that $B\cap A\cap E^*\cap C_\delta\notin I_\delta$.  We will exploit this as a way
of ensuring sets are $\id_p(\bar{C},\bar{I})$-positive.

By recursion on $i<\tau$, we will choose objects $B_i$, $\delta_i$, and $Y_i$ such that
\begin{enumerate}
\item $B_i\in M\cap\mathcal{P}(A)$
\sk
\item $\delta_i\in S\cap E^*$ (but not necessarily an element of $M$)
\sk
\item $|Y_i|$ is an $I_{\delta_i}$-positive subset of $B_i\cap E^*\cap C_{\delta_i}$ of cardinality $\cf(\mu)$
\sk
\item $B_i\cap\bigcup_{j<i}Y_j=\emptyset$.
\sk
\end{enumerate}

We postpone for now the proof that such objects can be found, and instead prove that the above construction
yields a collection $\langle D_i:i<\tau\rangle$ of disjoint $\id_p(\bar{C},\bar{I})$-positive sets. We do this by setting
\begin{equation*}
D_i:= B_i\setminus\bigcup_{i<k<\tau}B_k.
\end{equation*}
It should be clear that the collection $\{D_i:i<\tau\}$ is pairwise disjoint, so we are done if
we can prove that each $D_i$ is $\id_p(\bar{C},\bar{I})$-positive.

Note that since $M$ is closed under sequences of length $\tau$, the sequence $\langle B_i:i<\tau\rangle$
is an element of $M$, and therefore each $D_i$ is in $M$ as well. Requirement~(4) of
our construction tells us that $Y_i\subs D_i$. By~(3), there is an ordinal $\delta\in S$ (namely
$\delta_i$) for which
\begin{equation*}
D_i\cap A\cap E^*\cap C_\delta\notin I_\delta.
\end{equation*}
Since $D_i$ is in $M$, our ``litmus test'' tells us that $D_i$ must be $\id_p(\bar{C},\bar{I})$-positive, and so the collection
$\langle D_i:i<\tau\rangle$ contradicts the assumption that $\id_p(\bar{C},\bar{I})\cap \mathcal{P}(A)$ is weakly $\tau$-saturated. We are forced
to conclude that (\ref{413goal}) holds, and the theorem follows.

Now why can such objects be found?  Suppose that we have managed to find $B_j$, $Y_j$, and $\delta_j$ for
all $j<i<\tau$.  Since we assume (\ref{413goal}) fails and $|M|=(2^{\tau})^+$, we can find a collection $\langle X_\alpha:\alpha<(2^\tau)^+\rangle$
of sets in $M\cap\mathcal{P}(A)$ such that
\begin{equation*}
\alpha\neq\beta\Longrightarrow X_\alpha\neq_I X_\beta.
\end{equation*}
Since $i<\tau$ and  $Y_j$ is of cardinality at most $\cf(\mu)\leq\tau$ for each $j<i$, there exist distinct
$\alpha$ and $\beta$ such that
\begin{equation*}
X_\alpha\cap Y_j = X_\beta\cap Y_j\text{ for all }j<i.
\end{equation*}
Now either $X_\alpha\setminus X_\beta\notin I$ or $X_\beta\setminus X_\alpha\notin I$. In the former
case, let $B_i= X_\alpha\setminus X_\beta$, otherwise we set $B_i= X_\beta\setminus X_\alpha$.  In
either case,
\begin{equation}
\label{scriabin}
B_i\notin I,
\end{equation}
and
\begin{equation*}
B_i\cap Y_j=\emptyset\text{ for all }j<i.
\end{equation*}
Because of (\ref{scriabin}), we can find an ordinal $\delta_i\in S\cap E^*$ such that
\begin{equation}
\label{chopin}
B_i\cap E^*\cap C_{\delta_i}\notin I_{\delta_i}.
\end{equation}
We now choose $Y_i$ to be a subset of $B_i\cap E^*\cap\nacc( C_{\delta_i})$ that is cofinal in $\delta_i$
of order-type $\cf(\delta_i)=\cf(\mu)$ and such that the sequence $\langle \cf(\alpha):\alpha\in Y_i\rangle$
increases to $\mu$.  This can be done using (\ref{chopin}) and the definition of $I_{\delta_i}$, and with
this choice of $Y_i$, we have shown how to carry out the recursion.
\end{proof}

Returning to the matter of Theorem~\ref{reflthm}, we will need to make use of the following combinatorial notion
which appears in many guises in the literature.

\begin{defn}
Suppose $\theta$  is a regular cardinal, and $\mathcal{A}=\{A_\alpha:\alpha\in \Lambda\}$ is a family of sets,
 each of size $\theta$.  We say that $\mathcal{A}$ is {\em essentially disjoint} if there are sets $\{B_\alpha:\alpha\in\Lambda\}$
 such that each $B_\alpha$ is of cardinality less than $\theta$, and the family $\{A_\alpha\setminus B_\alpha:\alpha\in\Lambda\}$
 is pairwise disjoint.
\end{defn}

The following lemma is a sharpening of well-known results concerning indecomposable ultrafilters.

\begin{lem}
\label{essdisj}
Let $\theta<\kappa$ be regular cardinals, let  $I$ be an $\theta$-indecomposable ideal on~$\kappa$
extending the ideal of bounded sets, and let $\mathcal{A}$ be a family
of $\kappa$ sets each of cardinality $\theta$.  If every subfamily of $\mathcal{A}$ of size less than
$\kappa$ is essentially disjoint, then $\mathcal{A}$ can be written as a union of at most $|\mathcal{P}(\kappa)/ I|$
essentially disjoint families.
\end{lem}
\begin{proof}
Let $\mathcal{A}=\{A_\alpha:\alpha<\kappa\}$, and for each $\alpha<\kappa$ fix a function $F_\alpha$
witnessing that the family $\{A_\beta:\beta<\alpha\}$ is essentially disjoint, that is,
\begin{itemize}
\item $\dom(F_\alpha)=\alpha$,
\item $F_\alpha(\beta)\in [A_\beta]^{<\theta}$ for each $\beta<\alpha$, and
\item the family $\{A_\beta\setminus F_\alpha(\beta):\beta<\alpha\}$ is disjoint.
\end{itemize}
Also fix, for each $\alpha$, a bijection $b_\alpha$ between $A_\alpha$ and $\theta$.

For each $\beta<\kappa$ and $\epsilon<\theta$, let us define
\begin{equation}
B(\beta,\epsilon):=\{\alpha<\kappa: b_{\beta}[F_\alpha(\beta)]\subs\epsilon\}.
\end{equation}
For each $\beta<\kappa$, the sequence $\langle B(\beta,\epsilon):\epsilon<\theta\rangle$ is increasing, and
clearly
\begin{equation}
\bigcup_{\epsilon<\theta} B(\beta,\epsilon) = \kappa\setminus(\beta+1)\notin I.
\end{equation}
Since $I$ is $\theta$-indecomposable,
it follows that for each $\beta<\kappa$, there is a value $\epsilon(\beta)<\theta$ such that
 $B(\beta,\epsilon(\beta))$ is not in $I$.

We now define an equivalence relation on $\kappa$ according to the rule
\begin{equation}
\beta\sim\gamma \Longleftrightarrow B(\beta,\epsilon(\beta)) = B(\gamma,\epsilon(\gamma))\mod I.
\end{equation}
If we let $\tau$ denote the cardinality of $\mathcal{P}(\kappa)/I$, then it is clear that the number of
 $\sim$-equivalence classes is at most $\tau$.

Now define a function $F$ with domain $\kappa$ by
 \begin{equation}
 \label{fdef}
 F(\alpha):= b_\alpha^{-1}[\epsilon(\alpha)].
 \end{equation}
It is clear that $F(\alpha)$ is a subset of $A_\alpha$ of cardinality less than $\theta$.  To finish, we verify
that $A_\beta\setminus F(\beta)$ and $A_\gamma\setminus F(\gamma)$ are disjoint whenever $\beta\sim\gamma$.

Given $\beta\sim\gamma$, we note that
\begin{equation}
\label{equiv}
B(\beta,\epsilon(\beta))\cap B(\gamma,\epsilon(\gamma))\notin I.
\end{equation}
This is easily seen --- since the two sets are equivalent modulo $I$, the only way (\ref{equiv}) can fail
is if both are in $I$, but this would contradict the definition of $\epsilon(\beta)$.

In particular, this means that the two sets have non-empty intersection and so we can choose
 $\alpha\in B(\beta,\epsilon(\beta))\cap B(\gamma,\epsilon(\gamma))$.  Recall that the function $F_\alpha$
 has the property that
 \begin{equation}
\bigr(A_\beta\setminus F_\alpha(\beta)\bigl)\cap \bigr(A_\gamma\setminus F_\alpha(\gamma)\bigl)=\emptyset.
\end{equation}
We can appeal to (\ref{fdef}) and the definition of $B(\beta,\epsilon(\beta))$ to conclude that
$F_\alpha(\beta)$ is a subset of $F(\beta)$, and the same argument tells us that $F_\alpha(\gamma)\subs F(\gamma)$.
Therefore, the sets $A_\beta\setminus F(\beta)$ and $A_\gamma\setminus F(\gamma)$ are disjoint, as required.
\end{proof}

Now at last we are in a position to combine the preceding lemmas to give a proof
of Theorem~\ref{reflthm}.

\begin{proof}[Proof of Theorem~\ref{reflthm}]
Let $\lambda$, $\mu$, $S$, and $(\bar{C},\bar{I})$ be as in the statement of the theorem.  If $\id_p(\bar{C},\bar{I})$
is weakly $\mu$-saturated, then by Lemma~\ref{cardlem} we know that there is an $\id_p(\bar{C},\bar{I})$-positive
set $A$ for which
\begin{equation*}
|\mathcal{P}(A)/\id_p(\bar{C},\bar{I})|<\mu.
\end{equation*}

Suppose now that $T$ is a non-reflecting stationary subset of $\{\delta<\lambda:\cf(\delta)=\theta\}$ for some regular $\theta<\mu$
different from $\cf(\mu)$.  For each $\delta\in T$, let $A_\delta\subs\delta$ be cofinal of order-type $\theta$.
Since $T$ does not reflect, a well-known result (see Section~2.2 of~\cite{myhandbook}, for example) tells us that
every subset of $\{A_\delta:\delta\in S\}$ of cardinality less than $\lambda$ is essentially disjoint.

The ideal $\id_p(\bar{C},\bar{I})$ is $\theta$-indecomposable by Proposition~\ref{observation}, and it follows immediately
that the ideal $\id_p(\bar{C},\bar{I})\cap\mathcal{P}(A)$ is $\theta$-indecomposable as well.
 By an appeal to Lemma~\ref{essdisj}
we conclude that the family $\{A_\delta:\delta\in T\}$ is the union of at most $|\mathcal{P}(A)/\id_p(\bar{C},\bar{I})|$
essentially disjoint families.

Now  $|\mathcal{P}(A)/\id_p(\bar{C},\bar{I})|<\mu$, and therefore there must be a stationary $T^*\subs T$ for which
$\{A_\delta:\delta\in T^*\}$ is essentially disjoint.  This is absurd, as we immediately get a contradiction
to Fodor's lemma. Thus, the stationary set $T$ must reflect.
\end{proof}

\begin{cor}
Let $\lambda=\mu^+$ with $\mu$ a strong limit singular. If $\square_\mu$ holds, then for any
stationary $S\subs\{\delta<\lambda:\cf(\delta)=\cf(\mu)\}$ and $S$-fair pair $(\bar{C},\bar{I})$, we can
partition $\lambda$ into $\mu$ disjoint $\id_p(\bar{C},\bar{I})$-positive sets.
\end{cor}
\begin{proof}
The proof consists of the conjunction of Theorem~\ref{reflthm} with the well-known fact that $\square_\mu$ implies that every stationary subset
of $\lambda$ has a non-reflecting stationary subset.
\end{proof}

\section{Final Comments}

There are several natural questions raised by this research. One such is the question of the extent to which
similar results hold at successors of singular cardinals of countable cofinality, and this has been addressed in
recent joint work of the author and Saharon Shelah~\cite{819}.  In particular, we obtain coloring theorems (not quite as
strong as $\pr_1$) using $S$-fair pairs satisfying certain structural requirements.   These
$S$-fair pairs can be shown to exist at all successors of singular cardinals, including successors
of singular cardinals of countable cofinality. However, the following question is still very much open:

\begin{question}
Suppose $\lambda=\mu^+$ for $\mu$ singular of countable cofinality, and let $S\subs\{\delta<\lambda:\cf(\delta)=\cf(\mu)\}$.
Does there exist an $S$-good pair?
\end{question}
Of course, the main question we have not been able to answer is the following:
\begin{question}
Is it consistent that an ideal of the form $\id_p(\bar{C},\bar{I})$ (with $(\bar{C},\bar{I})$ an $S$-fair or even $S$-good pair)
can be weakly $\mu$-saturated?
\end{question}
The answer is almost certainly ``yes'', but we have made no inroads.
Finally, we have the following two questions. We will discuss the motivation in a moment.

\begin{question}
Suppose $\lambda=\mu^+$ for $\mu$ singular, and let $(\vec{\mu},\vec{f})$ be a scale for $\mu$, where
$\vec{\mu}=\langle \mu_i:i<\cf(\mu)\rangle$.  Define
\begin{equation*}
S^*=\{\alpha<\lambda:\cf(\alpha)=\mu_i\text{ for some }i<\cf(\mu)\}.
\end{equation*}
Can we find a stationary $S\subs\{\delta<\lambda:\cf(\delta)=\cf(\mu)\}$ and an $S$-fair pair $(\bar{C},\bar{I})$
such that $S^*\notin\id_p(\bar{C},\bar{I})$?
\end{question}

\begin{question}
Suppose $\lambda=\mu^+$ for $\mu$ singular, and let $(\bar{C},\bar{I})$ be an $S$-fair pair for some
stationary $S\subs\{\delta<\lambda:\cf(\delta)=\cf(\mu)\}$.  If there is a strictly increasing sequence
of regular cardinals $\vec{\mu}=\langle \mu_i:i<\cf(\mu)\rangle$ cofinal in $\mu$ such that
\begin{equation*}
\{\alpha<\lambda:\cf(\alpha)=\mu_i\text{ for some }i<\cf(\mu)\}\notin \id_p(\bar{C},\bar{I}),
\end{equation*}
must it be the case that this sequence $\vec{\mu}$ carries a scale for $\mu$?
\end{question}

The above questions are rooted in the observation that there are many examples in this area where the same
conclusions follow from hypotheses on club-guessing, as well as hypotheses concerning scales. For example, consider
the following pair of results\footnote{The first of these two results
appears explicitly as Conclusion~4.6 on page~73 of \cite{cardarith}, while the second follows easily from
the proof of Lemma~1.9 on page 121 of the same book.} due to Shelah:

Assume $(\vec{\mu},\vec{f})$ is a scale for $\mu$ and
\begin{equation*}
(\forall^*i<\cf(\mu))(\mu_i\nrightarrow [\mu_i]^{<\omega}_{\mu_i}).
\end{equation*}
Then we can conclude that $\mu^+\nrightarrow [\mu^+]^{<\omega}_{\mu^+}$ holds as well.

On the other hand, if we assume $(\bar{C},\bar{I})$ is an $S$-fair pair, and
\begin{equation*}
\{\alpha<\lambda: \cf(\alpha)\nrightarrow [\cf(\alpha)]^{<\omega}_{\cf(\alpha)}\}\notin\id_p(\bar{C},\bar{I}),
\end{equation*}
then we obtain the same conclusion $\mu^+\nrightarrow [\mu^+]^{<\omega}_{\mu^+}$.

Our questions are intended to probe the extent to which club-guessing and scales are related, and to see if
there is any deep reason why the we often have such pairs of results.  A positive answer to either Question~3
or Question~4 would be surprising as there don't seem to be any compelling reasons for the two ideas to
be connected. Still, the phenomena of ``same conclusion from parallel hypotheses'' is puzzling and perhaps there
is an explanation for it.


\begin{thebibliography}{10}



\bibitem{jb}
James E. Baumgarter.
\newblock A new class of order-types.
\newblock {\em Ann. Pure Appl. Logic}, 54(3):195--227, 1991.

\bibitem{bekkali}
Mohamed Bekkali.
\newblock {\em Topics in Set Theory}.
\newblock Number 1476 in Lecture Notes in Mathematics. Springer-Verlag, Berlin,
  1991.


\bibitem{cummings}
James Cummings.
\newblock Notes on singular cardinal combinatorics.
\newblock {\em Notre Dame J. Formal Logic}, 46(3):251--282, 2005.


\bibitem{535}
T.~Eisworth and S.~Shelah.
\newblock Successors of singular cardinals and coloring theorems {I}.
\newblock {\em Arch. Math. Logic}, 44(5):597--618, 2005.

\bibitem{819}
T.~Eisworth and S.~Shelah.
\newblock Successors of singular cardinals and coloring theorems {II}
\newblock { (submitted to {\em J. Symbolic Logic})}.


\bibitem{myhandbook}
Todd Eisworth
\newblock Successors of singular cardinals.
\newblock (chapter in the forthcoming {\em Handbook of Set Theory}).

\bibitem{ehr}
P.~Erd{\H{o}}s, A.~Hajnal, and R.~Rado.
\newblock Partition relations for cardinal numbers.
\newblock {\em Acta Math. Acad. Sci. Hungar.}, 16:93--196, 1965.


\bibitem{hsw}
M. Holz, K. Steffens, and E. Weitz.
\newblock{\em Introduction to cardinal arithmetic}
\newblock Birkh\"{a}user Advanced Texts: Basler Lehrb\"{u}cher, Birkh\"{a}user Verlag, Basel, 1999.


\bibitem{cardarith}
Saharon Shelah.
\newblock {\em {Cardinal Arithmetic}}, volume~29 of {\em {Oxford Logic
  Guides}}. {Oxford University Press}, 1994.

\bibitem{Sh:355}
Saharon Shelah.
\newblock {$\aleph _{\omega +1}$ has a Jonsson Algebra}.
\newblock In {\em {Cardinal Arithmetic}}, volume~29 of {\em {Oxford Logic
  Guides}}, Chapter~II. {Oxford University Press}, 1994.

\bibitem{Sh:365}
Saharon Shelah.
\newblock {There are Jonsson algebras in many inaccessible cardinals}.
\newblock In {\em {Cardinal Arithmetic}}, volume~29 of {\em {Oxford Logic
  Guides}}, Chapter III. {Oxford University Press}, 1994.

\bibitem{572}
Saharon Shelah.
\newblock Colouring and non-productivity of {$\aleph\sb 2$}-c.c.
\newblock {\em Ann. Pure Appl. Logic}, 84(2):153--174, 1997.

\bibitem{413}
Saharon Shelah.
\newblock More {J}onsson algebras.
\newblock {\em Arch. Math. Logic}, 42(1):1--44, 2003.

\bibitem{stevochapter}
Stevo Todor{\v{c}}evi{\'c}.
\newblock Coherent sequences.
\newblock Chapter in the forthcoming Handbook of Set Theory.

\bibitem{minimal}
Stevo Todor{\v{c}}evi{\'c}.
\newblock Partitioning pairs of countable ordinals.
\newblock {\em Acta Math.}, 159(3-4):261--294, 1987.

\bibitem{stevobook}
Stevo Todor{\v{c}}evi{\'c}.
\newblock {\em {Walks on ordinals and their characteristics}}, volume~263 of {\em {Progress in Mathematics}}. {Birkhauser}, 2007.

\end{thebibliography}
\end{document}